\documentclass{amsart}

\usepackage{amssymb} 
\usepackage[mathscr]{eucal}
\usepackage{graphicx}

\newtheorem{theorem}{Theorem}
\newtheorem{lemma}[theorem]{Lemma}
\newtheorem{proposition}[theorem]{Proposition}
\newtheorem{corollary}[theorem]{Corollary}

\def\half{{\textstyle{\frac12}}}

\def\fourth{{\textstyle{\frac14}}}
\def\onesixteenth{{\textstyle{\frac1{16}}}}
\def\oneoverF{\textstyle{\frac{1}{F}}}
\def\pa{\partial}
\def\hG{ \mathcal{G} }

\def\P{\mathscr{P}}
\def\drho{\dot\rho}
\def\dP{\dot{\mathscr{P}}}
\def\ddrho{\ddot\rho}
\def\ddP{\ddot{\mathscr{P}}}
\def\r0{\rho_o}
\def\dr0{\dot{\rho_o}}
\def\ddr0{\ddot{\rho_o}}
\def\R{\mathbb{R}}
\def\I{(-\varepsilon,\varepsilon)}
\def\L{\mathcal L}
\def\C{\mathcal C}
\def\T{\mathcal T}

\newcounter{cl}
\newenvironment{clist}[1]
  {
    \begin{list}{#1}
      {\usecounter{cl}
       \setlength{\topsep}{0pt}
       \setlength{\labelsep}{4pt}
       \setlength{\labelwidth}{10pt}
       \setlength{\leftmargin}{15pt}
      }
  }
  {
    \end{list}
  }

\begin{document}


\title[Randers geodesics]
      {Geodesics in Randers spaces \\ of constant curvature} 

\author[Robles]
       {Colleen Robles}

\address{Department of Mathematics, \\
         University of Rochester    \\
         Rochester, NY 14607, USA}

\email{robles@math.rochester.edu}




\begin{abstract}
       Geodesics in Randers spaces of constant curvature are classified.
\end{abstract}


\maketitle 


\section{Introduction} 

  Randers metrics have received much attention lately as solutions to 
  Zermelo's problem of navigation; largely because this navigation structure 
  provides the framework for a complete classification of constant flag 
  curvature Randers spaces.  ({\it Flag curvature} is the Finslerian analog of 
  Riemannian sectional curvature.  See \cite{BR04}.)  Briefly, a Randers 
  metric is of constant flag curvature if and only if it solves Zermelo's 
  problem of navigation on a Riemannian manifold of constant sectional 
  curvature under the influence of an infinitesimal homothety $W$.  See 
  Subsection \ref{sec:randers} for a sketch of the navigation problem, and 
  Theorem \ref{thm:class} for an explicit 
  statement of the classification result.

  The aim of this paper is to develop a geometric description of the geodesics
  in these spaces of constant curvature.  Intuitively, these paths minimize 
  travel time across a Riemannian landscape under windy conditions.  Presently
  I will show that these curves are given by composing geodesics of the 
  Riemannian metric with the flow generated by $W$.  This claim is 
  formalized by Theorem \ref{thm:geodesics}.  

  Geodesics on surfaces of constant, nonpositive curvature are 
  illustrated in Section \ref{sec:pictures}.  We then turn, in 
  \S\ref{sec:sphere}, to the constant flag curvature $K=1$ Randers metrics on 
  $S^n$.  The case of the sphere is especially interesting; it is possible to 
  endow this closed manifold with a metric whose geodesics display distinctly 
  non-Riemannian behaviors.  For example,
  \begin{clist}{(\arabic{cl})}
    \item
      A metric is {\it projectively flat} if every point admits coordinates 
      in which the geodesics are straight lines.  Beltrami's theorem assures 
      us that {\it a Riemannian metric is of constant sectional curvature if 
      and only if it is projectively flat}.  In contrast few Randers spaces of 
      constant flag curvature are projectively flat.  There are infinitely many
      non-isometric Randers metrics of constant, {\it positive} curvature.  
      Only the trivial (Riemannian) metrics are projectively flat.  See Section
      7 of \cite{BRS04} for a thorough discussion. 
    \item 
      On the Riemannian $n$-sphere all geodesics close with length 
      $2 \pi$.  In contradistinction, it is possible to equip the sphere with 
      a non-Riemannian Randers metric of constant curvature $K=1$ for which 
      either (i) all the geodesics close, or (ii) only finitely many of the 
      geodesics close!  In either case, the geodesics will not all have the 
      same length.  This is a consequence of the fact that non-Riemannian 
      Randers metrics are not {\it reversible} (a.k.a. {\it symmetric}):  in 
      general, the Randers length $F(V)$ of a vector $V$ is not equal to the 
      length $F(-V)$ of $-V$.
  \end{clist}
\vspace{0.1in}

  Studies of conjugate and cut points occupy the last two sections of the 
  paper.  Section \ref{sec:conjugate_points} establishes a correspondence 
  between points conjugate to $p$ with respect to the Riemannian metric, and 
  points conjugate to $p$ with respect to the Randers metric solving the 
  navigation problem on the Riemannian manifold.  The details are spelled out 
  in Proposition \ref{prop:conjugate_pts}.

  The final section of the paper investigates minimizing geodesics.  We will 
  see that a Randers geodesic is a global minimizer if and only if its 
  associated Riemannian geodesic is as well.  See Proposition 
  \ref{prop:global_minimizers} for a precise statement.

  The remainder of the introduction is given to a brief discussion of Randers 
  metrics as solutions to Zermelo's problem of navigation.  This material has 
  been carefully presented elsewhere (references are provided), so I will 
  present only a cursory introduction with the principle goal of establishing 
  notation.  To begin,
  \begin{itemize}
    \item[$\circ$] Points on the $n$-dimensional manifold $M$ are denoted by 
          $p$ or $x$.
    \item[$\circ$] Tangent vectors are given by $y \in T_p M$, with components 
          $y^i \partial_{x^i}$ relative to local coordinates $x=(x^i)$ on $M$.
    \item[$\circ$] Partial derivatives are denoted by the subscripts ${x^i}$ 
          and ${y^i}$.
  \end{itemize}


  \subsection{Randers metrics and Zermelo Navigation}
  \label{sec:randers}

  In 1941 G. Randers \cite{Ra41} introduced a Finsler metric by modifying 
  a Riemannian metric $a := a_{ij}(x) \, dx^i \otimes dx^j$ by a linear term
  $b := b_{i}(x) \, dx^i$.  The resulting Minkowski norm on $T_xM$ is given by 
  \begin{displaymath}
    F(x,y) \ := \ \alpha(x,y) \, + \, \beta(x,y) \ 
              = \ \sqrt{a_{ij}(x) y^i y^j} \, + \, b_i(x) y^i \, , 
    \quad y = y^i \pa_{x^i} \in T_x M \, . 
  \end{displaymath}
  By requiring $a(b,b) < 1$, we ensure that $F$ is positive.  This simple 
  condition also guarantees that the metric is strongly convex.  That is, the 
  Hessian $g_{ij}(x,y) := (\half F^2)_{y^i y^j}$ is positive definite for all 
  nonzero $y$.  See \cite{BCS00, BR04} for this result and a through 
  treatment of Randers metrics.

  Z. Shen \cite{S02} has identified Randers metrics with solutions to a 
  navigation problem.  In 1931, Zermelo posed and answered the following 
  question \cite{Z31,C99}:  Suppose a ship sails the sea on calm waters.  
  Imagine a mild breeze comes up.  How must the captain guide the ship to 
  reach a given destination in the shortest time?  

  Zermelo assumes that the sea is $\R^2$, with the flat/Euclidean metric.
  Shen extended the seascape to an arbitrary Riemannian manifold 
  $(\mathcal M,h)$.  Assuming a time-independent wind, he found that the paths 
  minimizing travel-time are exactly the geodesics of a Randers metric
  \begin{equation}
  \label{eqn:F}
    F(x,y) = \alpha(x,y) + \beta(x,y)
    = \frac{\sqrt{\lambda |y|^2 + W_0{}^2}}{\lambda}
      - \frac{W_0}{\lambda} \, .
  \end{equation}
  Here $W = W^i \pa_{x^i}$ is the velocity vector field of the wind, 
  \begin{center}
  $|y|^2 = h(y,y)$ , $\quad \lambda = 1 - |W|^2 \quad$ and 
  $\quad W_0 = h(W,y)$.
  \end{center}  
  We say {\it $F$ solves Zermelo's problem of navigation}.  The defining 
  Riemannian metric $a$ and 1-form $b$ of the Randers metric are 
  \begin{displaymath}
    a_{ij} = \frac{\lambda h_{ij} + W_i W_j}{\lambda^2} 
    \quad \hbox{and} \quad
    b_i = -\frac{W_i}{\lambda} \, , \quad \hbox{where } W_i := h_{ij} W^j \, .
  \end{displaymath} 
  Requiring $|W| < 1$ ensures that this 
  $a_{ij}$ is indeed a positive definite Riemannian metric.  Additionally, 
  $h(W,W) = a(b,b)$.  Therefore, the condition $|W|<1$ also ensures the 
  positivity of $F$.  Consequently $F$ is defined on the open sub-manifold 
  $\{ |W| < 1 \} \subset \mathcal M$.  
  A straightforward computation establishes
  \begin{lemma}
  \label{lem:length}
    $F(y) = 1$ $\Longleftrightarrow$ $|y-W| = 1$.  Geometrically, this means
    the unit sphere of $F$ in $T_xM$ differs from the unit sphere of $h$ 
    by a translation along $W(x)$.
  \end{lemma}

  At this point it is natural to ask if {\it every} Randers space $(M,F)$ 
  arises as the solution to Zermelo's problem of navigation for a canonical
  choice of $(h,W)$.  The answer is yes \cite{BR04}.  As a result, 
  {\it Randers metrics are naturally identified with solutions to the 
  navigation problem.}


  \subsection{The result}
  \label{sec:result}

  It is now possible to state the main result of the paper.  The notation
  $\mathcal L$ below denotes Lie differentiation.
  \begin{theorem}
  \label{thm:geodesics}
    Assume $(\mathcal M , h)$ is a Riemannian manifold equipped with an 
    infinitesimal homothety $W$;  $\mathcal L_W h = \sigma h$, $\sigma$ a 
    constant.  Let $F$ denote the Randers metric solving Zermelo's problem of 
    navigation on $M := \{ |W| < 1 \} \subset \mathcal M$.  Then the unit 
    speed geodesics $\P : (-\varepsilon , \varepsilon) \to M$ of $F$ are given 
    by $\P(t) = \varphi(t,\rho(t))$.  Here
    \begin{itemize}
      \item[$\circ$]
        $\rho : (-\varepsilon , \varepsilon) \to \mathcal M$ is a geodesic of 
        $h$ parameterized so that $|\drho(t)|^2 = e^{-\sigma t}$,
      \item[$\circ$]
        shrinking $\varepsilon$ if necessary, 
        $\varphi : (-\varepsilon, \varepsilon) \times U \to M$ is the 
        flow of $W$ defined on a neighborhood $U$ of $\rho(0)$ so that 
        $\rho(t) \in U$, for all $t \in (-\varepsilon, \varepsilon)$.
    \end{itemize}
  \end{theorem}
  In the abstract I promised to classify the geodesics in Randers spaces of 
  constant flag curvature.  The theorem does just that.  This 
  is because a Randers metric has constant flag curvature if and only if 
  it solves the navigation problem on a Riemannian space of constant 
  curvature under the influence of an infinitesimal homothety, and the 
  geodesic structure of Riemannian space forms is standard knowledge.
  \begin{theorem}[Constant flag curvature classification \cite{BRS04}]
  \label{thm:class}
    Let $F$ be a Randers metric on a manifold $M$ with navigation data 
    $(h,W)$, $|W|<1$.  Then $F$ is of constant flag curvature $K$ if and only 
    if 
    \begin{itemize}
      \item[(1)]
        The Riemannian space $(M,h)$ is of constant sectional curvature 
        $K + \onesixteenth \sigma^2$.
      \item[(2)] 
        The vector field $W$ is an infinitesimal homothety of $h$, 
        $\mathcal L_W h = \sigma h$.  
    \end{itemize}
    The constant $\sigma$ vanishes when $h$ is not flat, and $W$ is an 
    infinitesimal isometry.
  \end{theorem}

  Note that Theorem \ref{thm:geodesics} allows an arbitrary Riemannian metric 
  $h$.  As a result it describes the geodesics of a larger class of Randers 
  metrics than simply those of constant flag curvature.  In this context 
  the theorem may only be considered a `description,' not a `classification,'
  since the geodesic structure of these more general Riemannian metrics is 
  typically not known.


  \subsection{Geodesic equations}
  \label{sec:geodesics}
  A curve $\rho : \I \to \mathcal M$ is a geodesic of $h$ if it satisfies the 
  {\it geodesic equation}
  \begin{equation}
  \label{eqn:hgeo}
    \ddrho^i + 2 \, \hG^i(\rho,\drho) 
    = \frac{d}{dt}\left( \ln | \drho | \right) \drho^i \, .
  \end{equation}
  Note that the right-hand side is zero when $\rho$ is parameterized with 
  constant speed.  The geodesic spray coefficients $\hG^i$ are given by 
  \begin{displaymath}
    \hG^i(x,y) = \half \gamma^i{}_{jk}(x) y^j y^k \, ,
  \end{displaymath}
  where 
  $\gamma^i{}_{jk} = \half h^{is} ( h_{sj,x^k} - h_{jk,x^s} + h_{ks,x^j} )$
  denotes the Christoffel symbols of $h$.

  Similarly, a curve $\P : \I \to M$ will be a geodesic of the Randers metric 
  $F$ if it satisfies the {\it geodesic equation}
  \begin{displaymath}
    \ddP^i + 2 G^i(\P,\dP) = \frac{d}{dt}\left( \ln F(\dP) \right) \dP^i \, .
  \end{displaymath}
  (See \cite{BCS00} for a discussion of Finslerian geodesics.)
  The geodesic coefficients of $F$ are related to those of 
  the Riemannian metric $a$ by equation (11.3.12) of \cite{BCS00}.  
  ({\it Beware!}  Their $G$ is {\it twice} ours.).  Similarly, the 
  relationship between the geodesic spray coefficients of 
  $a$ and $h$ is derived on page 235 of \cite{BR04}.  As a result, the 
  geodesic coefficients of $F$ are related to the those of $h$ by
  \begin{displaymath}
    G^i = \hG^i + \zeta^i \, ,
  \end{displaymath}
  where
  \begin{equation}
  \label{eqn:zeta}
    \zeta^i = \fourth \left( \oneoverF y^i - W^i \right) 
                      \left( 2 F \mathcal S_0 - \mathcal L_{00}
                             - F^2 \mathcal L_{WW} \right) 
            - \fourth F^2 \left( \mathcal S^i + \mathcal T^i \right)
            - \half F \mathcal C^i{}_0 \, .
  \end{equation}
  Above I have utilized the notation
  \begin{displaymath}
  \begin{array}{rclcrcl}
    \L_{ij} & = & W_{i:j} + W_{j:i}
    & \quad &
    \C_{ij} & = & W_{i:j} - W_{j:i} \\
    \mathcal S_i & = & W^s \L_{si}
    & \quad &
    \T_i & = & W^s \C_{si} \, .
  \end{array}
  \end{displaymath}
  The colon `:' denotes covariant differentiation so that 
  $W_{i:j} = W_{i,x^j} - W_s \gamma^s{}_{ij}$.
  Indices on these tensors are raised with the inverse of $h$.  For example,
  $\mathcal S^i = h^{ij} \mathcal S_j$.  As before the subscript 0 denotes 
  contraction with $y$, $\C^i{}_0 = h^{ij} \C_{jk} y^k$.  Finally, 
  $\L_{WW} = W^i W^j \L_{ij}$.


\section{Proof of Theorem \ref{thm:geodesics}}
\label{sec:proof}

  \subsection{Preliminaries}
  \label{sec:prelim}

  We are given three items.  First, a Randers space $(M , F)$ induced by 
  Zermelo navigation on the Riemannian manifold $(\mathcal M,h)$ under the 
  influence of an infinitesimal homothety $W$;
      \begin{equation}
      \label{eqn:lie}
        \L_W h = \sigma h \, , \quad \hbox{ $\sigma$ a constant.}
      \end{equation}
  Second, a geodesic of $h$, $\rho : \I \to U \subset \mathcal M$, 
  parameterized so that $|\drho(t)|^2 = e^{-\sigma t}$.  And third, the flow 
  $\varphi : \I \times U \to M$ of $W$, defined on a neighborhood $U$ of 
  $\rho(0)$.  There is no loss of generality in assuming that $U$ admits 
  coordinates $(x^i)$.

  With these data in hand define a curve $\P : \I \to M$ by 
  \begin{displaymath}
    \P(t) = \varphi(t,\rho(t)) \, .
  \end{displaymath}
  The proof is complete once we establish 

  (i) that the curve $\P$ is parameterized with unit $F$-speed, $F(\dP) = 1$; 
  and 
  
  (ii) $\P$ is a geodesic of $F$, 
  $\ddP + 2 G(\P,\dP) \, = \, \ddP + 2 \, \hG(\P,\dP) + 2 \zeta(\P,\dP) 
   = \, 0$.

  \noindent  
  Begin by computing 
  \begin{displaymath}
    \dP = W + d \varphi ( \drho ) \, ,
  \end{displaymath}
  where $d \varphi$ denotes the differential of the map 
  $\varphi( t , \cdot ) : U \to M$.  With respect to local 
  coordinates $(x^i)$ this reads $\dP^i = W^i + \varphi^i{}_{,x^j}\, \drho^j$.
  Similarly, 
  \begin{equation}
  \label{eqn:ddP}
    \ddP^i = W^i{}_{,x^j} W^j 
           + 2 W^i{}_{,x^j} \, \varphi^j{}_{,x^k} \, \drho^k 
           + \varphi^i{}_{,x^j x^k} \, \drho^j \drho^k 
           + \varphi^i{}_{,x^j} \, \ddrho^j \, .
  \end{equation}


  \subsection{Implications of the homothety hypothesis}
  \label{sec:implications}

  First, the differential $d\varphi$ of 
  $\varphi(t,\cdot) : U \to M$ scales vectors by a conformal factor of 
  $e^{\sigma t/2}$:  $| d \varphi \, V |^2 = e^{\sigma t} |V|^2$.
  In particular, 
  $| d \varphi ( \drho ) |^2 = e^{\sigma t} | \drho(t) |^2 = 1$.
  Lemma \ref{lem:length} implies $F(\dP) = 1$, and (i) is established.

  Second, since $\varphi( t_o , \cdot )$ is a homothety, 
  $\r0(t) := \varphi(t_o,\rho(t))$ is an geodesic of $h$ with speed 
  $|\dr0| = e^{\sigma( t_o - t )/2}$.  The geodesic equation (\ref{eqn:hgeo})
  implies
  \begin{equation}
  \label{eqn:homhyp}
    \qquad \qquad
    \ddr0^i + 2 \, \hG^i(\r0,\dr0) = \frac{d}{dt} 
                                     \left[ \half \sigma( t_o - t ) \right]
                                     \dr0^i \, .
  \end{equation}
  Differentiating by $t$ yields
    $\dr0^i = \varphi^i{}_{,x^j} \, \drho^j$ and  
    $\ddr0^i = \varphi^i{}_{,x^j x^k} \, \drho^j \drho^k
            + \varphi^i{}_{x^j} \, \ddrho^j$,
  where the partial derivatives of $\varphi$ are evaluated at $(t_o, \rho(t))$.
  Observe that $\r0(t_o) = \P(t_o)$ and $\dr0(t_o) = \dP(t_o) - W_{\P(t_o)}$.
  Consequently, at $t = t_o$, Equation \ref{eqn:homhyp} reads
  \begin{displaymath}
      \varphi^i{}_{,x^j x^k} \, \drho^j \drho^k
    + \varphi^i{}_{,x^j} \, \ddrho^j 
    + 2 \, \hG^i(\P,\dP - W) = -\half \sigma ( \dP^i - W^i ) \, .
  \end{displaymath}
  Since $t_o$ is arbitrary, this expression must hold for all $t$.  Making use
  of the equality
  $\half( \mathcal L^i{}_{j} + \mathcal C^i{}_j ) = W^i{}_{:j} =
  W_{i,x^j} + W^i \, \gamma^i{}_{jk}$, and the fact that 
  $\hG^i(x,y) = \half \, \gamma^i{}_{jk}(x) y^j y^k$ is quadratic in $y$, 
  this updates the formula (\ref{eqn:ddP}) for $\ddP$ to 
  \begin{eqnarray*}
    \ddP^i & = & -2 \, \hG^i(\P,\dP) + \half \sigma ( W^i - \dP^i ) 
             + \mathcal L^i{}_0 + \mathcal C^i{}_0 
             - \half \mathcal S^i + \half \mathcal T^i \\
           & = & -2 \, \hG^i(\P,\dP) + \half \sigma \dP^i + \half \mathcal T^i 
             + \mathcal C^i{}_0 \, .
  \end{eqnarray*}
  In the second equality I have made use of the relationships 
  \begin{center}
        $\mathcal L^i{}_0 = \sigma \dP \quad $ 
    and $\quad \mathcal S^i = \sigma W^i$,
  \end{center}
  both consequences of the facts that $W$ is an infinitesimal homothety 
  (Equation \ref{eqn:lie}), and the components of the Lie derivative are 
  $(\L_Wh)_{ij} = \L_{ij}$.


  \subsection{Conclusion}

  Finally, we turn to the expression (\ref{eqn:zeta}) for $\zeta$ in 
  Subsection \ref{sec:geodesics}.
  Note that Equation \ref{eqn:lie} reduces the term 
  $(2 F \mathcal S_0 - \mathcal L_{00} - F^2 \mathcal L_{WW})$ to 
  $-\sigma |\dP - W|^2$, which is equal to $-\sigma$ by (i) and Lemma 
  \ref{lem:length}.  Hence, 
  \begin{displaymath}
    \zeta^i(\P,\dP) = - \fourth \sigma \dP^i - \fourth \mathcal T^i 
                      - \half \mathcal C^i{}_0 \, . 
  \end{displaymath}
  A quick comparison with the final form of $\ddP$ above yields 
  $\ddP = -2 \, \hG - 2 \zeta = - 2 G$.  Hence (ii) holds, and $\P$ is a 
  geodesic of $F$. \hfill Q.E.D.
  

  \subsection{The general case}

  Theorem \ref{thm:geodesics} does not hold for arbitrary $W$.  Consider, as a 
  counter-example, the vector field $W = \varepsilon \pa_\phi$ 
  ($0 < \varepsilon < 1$) defined on the 
  dense open subset of $S^2$ covered by the coordinate map $(\theta , \phi) 
  \mapsto ( \sin \theta \, \sin \phi , \cos \theta \, \sin \phi , \cos \phi )$.
  Note that $W$ is not an infinitesimal isometry of standard Riemannian metric 
  $h$ on the sphere, and there is no parameterization $\rho(t)$ of the equator 
  $\{ \phi = \frac{\pi}{2} \}$ (a geodesic of $h$) with the property that 
  $\P(t) = \varphi( t , \rho(t) )$ is a unit speed geodesic
  of the associated Randers metric.


\section{Photo gallery for nonpositive curvature}
\label{sec:pictures}

This section contains a selection of graphics illustrating solutions to 
Zermelo's problem of navigation (i.e. geodesics of $F$) on subsets of  
Euclidean plane and the Poincar\'e disc.  The subset $M$ on which $F$ is 
defined will be determined by the constraint $|W|<1$.  These metrics are 
of constant, non-positive curvature.


  \subsection{In the Euclidean plane}
  \label{sec:plane}

  We begin with the Euclidean plane $(\mathbb{R}^2 , h)$.  The infinitesimal 
  homotheties are  given by 
  \begin{displaymath}
    W(u,v) = \left(
               \begin{array}{c}
                 \half \sigma u + kv + c_1 \vspace{0.05in} \\
                 \half \sigma v - ku + c_2
               \end{array}
             \right) \, , \quad 
  (u,v)\in\mathbb{R}^2 \, .  
  \end{displaymath}
  Here $k$, $c_1$, $c_2 \in \mathbb{R}$ are constants.  Theorem \ref{thm:class}
  tells us the Randers metric $F$ solving Zermelo's problem of navigation for 
  $(h,W)$ is of constant flag curvature $K = -\onesixteenth \sigma^2 \le 0$.

  If $k = 0$, then $F$ is projectively flat (see Section 7 of \cite{BRS04}).  
  These spaces admit coordinate systems in which the geodesics are straight 
  lines.  In the examples below I will restrict attention to examples that are 
  not projectively flat, i.e. infinitesimal homotheties with $k \not= 0$.

  {\bf Case I: $W$ an infinitesimal rotation}.
  Set $\sigma = 0$.  Then the Randers metric has constant flag curvature $K=0$.
  Let's suppose $c_i = 0$ for simplicity.  The infinitesimal rotation
  $W(u,v) = (v,-u)$ has norm less than 1 on the unit disc 
  $M = \{ (u,v) \in \mathbb{R}^2 \ : \ u^2 + v^2 < 1 \}$.  

  The geodesic characterization of Theorem \ref{thm:geodesics} tells us to 
  parameterize the Euclidean line $\rho$ with unit speed.  
  Take $\rho(t) = ( t + u_o , v_o )$.  Since $\rho(-t)$ is also a unit 
  speed geodesic of $h$, both $\P_+ = \varphi(t,\rho(t))$ and 
  $\P_-(t) = \varphi(t, \rho(-t))$ are geodesics of the Randers metric $F$.

  \begin{center}
  \renewcommand{\arraystretch}{1.5}
  \begin{tabular}{cp{10mm}c}
    \includegraphics{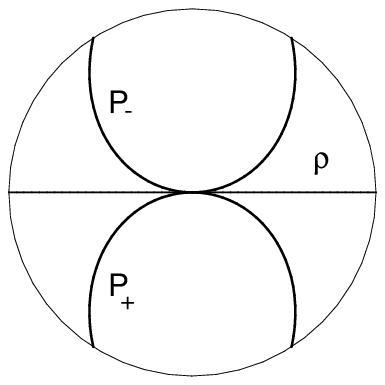} & &
    \includegraphics{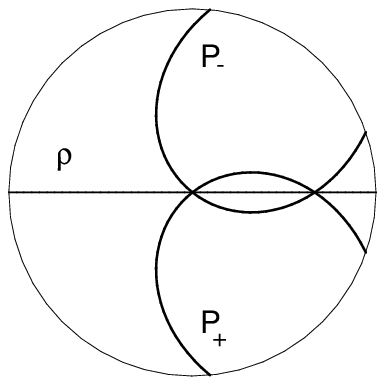} \\
    Example 1: $u_o = 0 = v_o$             & &
    Example 2: $u_o = \frac{2}{3}$, $v_o = 0$ \\
    \includegraphics{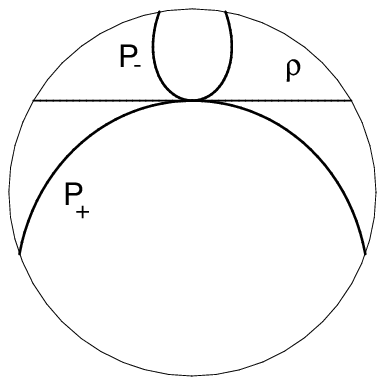} & &
    \includegraphics{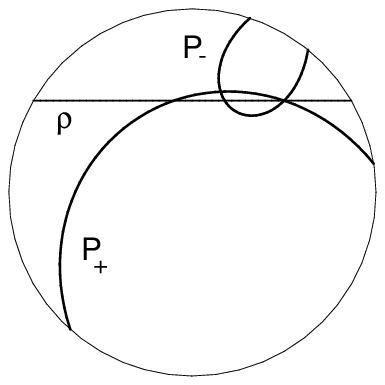} \\
    Example 3: $u_o = 0$, $v_o = \half$    & &
    Example 4: $u_o = \half = v_o$      
  \end{tabular}
  \end{center}
  Note that the two $F$-geodesics $\P_+$ and $\P_-$ in Example 1 are tangent 
  at the center of the disc.  Their initial vectors at this point are 
  $\dP_+(0) = (1,0)$ and $\dP_-(0) = (-1,0)$.  In Riemannian geometry two 
  geodesics passing through a common point in opposite directions necessarily 
  trace the same curve.  In fact, all reversible Finsler metrics have this 
  property.  As Examples 1 and 3 above indicate, the 
  phenomenon does not extend to the non-reversible setting.  In general, 
  if $\P(t)$ is a geodesic, $\P(-t)$ will not be a geodesic.

  {\bf Case II: {\boldmath$W$} an infinitesimal homothety}.
  \label{CaseII}  
  When $\sigma$ is nonzero the 
  metric $F$ has negative flag curvature $K = -\onesixteenth \sigma^2$.  
  In this case a suitably chosen translation effects a change of 
  coordinates for which $c_i = 0$.  The infinitesimal homotheties are of the 
  form $W(u,v) = -\half \sigma (u,v) + k (v , -u)$.    According to Theorem 
  \ref{thm:geodesics} we must parameterize the line $\rho$ so that 
  $|\rho(t)|^2 = e^{-\sigma t}$.  To that end take 
  $\rho_+(t) = [ \frac{2}{\sigma}( e^{-\sigma t / 2} - 1 ) + u_o , v_o ]$.
  The curve 
  $\rho_-(t) = [ \frac{2}{\sigma}( 1 - e^{-\sigma t / 2} ) + u_o , v_o ]$ 
  traverses the same path, again with Euclidean speed $e^{-\sigma t / 2}$, but
  in the opposite direction.  Hence, $\P_+(t) = \varphi(t,\rho_+(t))$ and 
  $\P_-(t) = \varphi(t , \rho_-(t))$ are $F$-geodesics, both generated from 
  the same Euclidean line.

  For $\sigma = \sqrt{2}$ and $k = 1/\sqrt{2}$  we have $|W|<1$ when 
  $u^2 + v^2 < 1$.  Hence $M \subset \mathbb{R}^2$ is again the unit disc.
  \begin{center}
  \begin{tabular}{cp{10mm}c}  
    \includegraphics{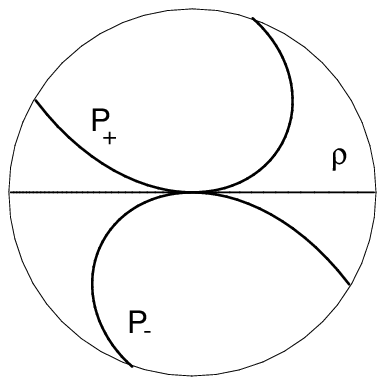} & &
    \includegraphics{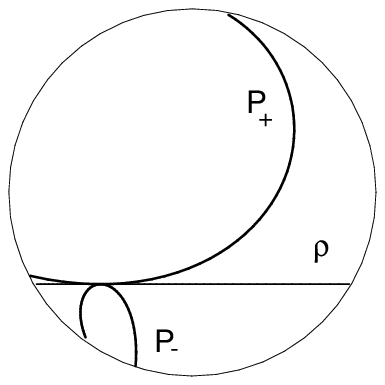} \\
    Example 5: $u_o = 0 = v_o$    & &
    Example 6: $u_o = -\half = v_o$ 
  \end{tabular}
  \end{center}


  \subsection{On the Poincar\'e disc}

  Next let's turn our attention to the Poincar\'e model $(D^2 , h)$
  of hyperbolic geometry on the unit disc with constant sectional curvature -1.
  In this case $\sigma$ must vanish; the infinitesimal homotheties $W$ are 
  necessarily infinitesimal isometries, and the resulting Randers metric $F$ 
  is of constant flag curvature $K=-1$ (Theorem \ref{thm:class}).  The metric 
  $F$ will be projectively flat if and only if $F$ is Riemannian ($W = 0$) 
  \cite{BRS04}.

  Since $\sigma = 0$ we parameterize the 
  geodesics of $h$, arcs of circles intersecting $\partial D$ 
  orthogonally, with unit hyperbolic speed.  If $\rho(t)$ is one such 
  parameterization, then both $\P_+(t) = \varphi(t,\rho(t))$ and 
  $\P_-(t) = \varphi(t,\rho(-t))$ are geodesics of $F$.  In Examples 7 and 9,
  $\rho(t) = ( 0 , \tanh(t/2) )$.  The geodesic $\rho$ of Examples 8 and 10 is 
  a hyperbolic translation of this vertical line $( 0 , \tanh(t/2) )$ along 
  the horizontal $u$-axis.

  {\bf Case III: {\boldmath$W$} is an infinitesimal rotation}.
  \begin{center}
  \begin{tabular}{cp{10mm}c}  
    \includegraphics{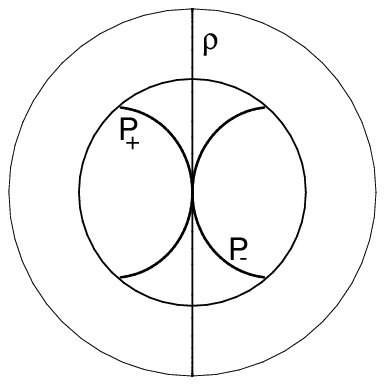} & &
    \includegraphics{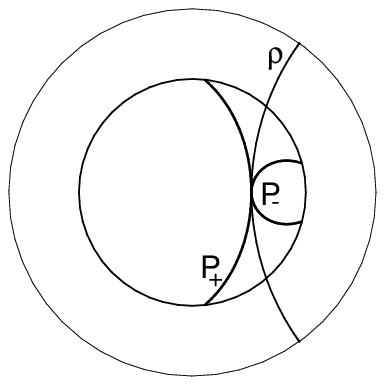} \\
    Example 7                           & &
    Example 8
  \end{tabular}
  \end{center}
  In the two examples above we take $W$ to be an infinitesimal rotation 
  about the origin, $W(u,v) = \half (v,-u)$.  The resulting Randers metric 
  is defined on $M$ the disc of radius $(\sqrt{5}-1)/2$, determined by the 
  condition $|W|<1$.  

  {\bf Case IV: {\boldmath$W$} is an infinitesimal hyperbolic translation}.
  Here the flow $\varphi(t, \cdot) : \mathbb{D} \to \mathbb{D}$ is a 
  hyperbolic translation fixing the horizontal $u$-axis.  The corresponding 
  vector field $W$ has norm less than one inside the smaller, eye-shaped 
  region $M$ indicated below.
  \begin{center}
  \begin{tabular}{cp{10mm}c}  
    \includegraphics{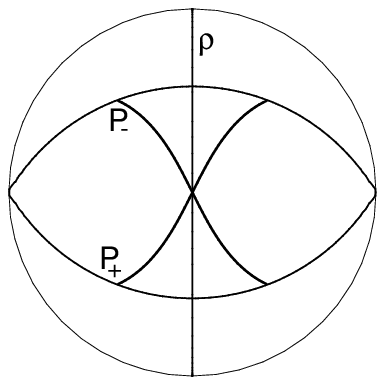} & &
    \includegraphics{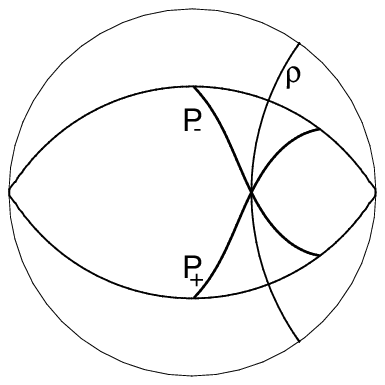} \\
    Example 9                           & &
    Example 10
  \end{tabular}
  \end{center}


  \subsection{When $(M , h)$ is complete}
  \label{sec:closed}  
  In each of the examples above $(M,h)$ fails to be complete.  
  \begin{proposition}
  Suppose is $(M,F)$ of constant nonpositive flag curvature, and 
  $(M,h)$ is complete.  Then $F$ is either locally Minkowski, or 
  a negatively curved Riemannian space form.
  \end{proposition}
  \begin{proof} 
  Theorem \ref{thm:class} assures us that the Riemannian space 
  $(M,h)$ must be of constant sectional curvature $\kappa \le 0$.  Since 
  $(M,h)$ is complete it admits a universal cover $(\tilde M , \tilde h)$ 
  by Euclidean $n$-space $\mathbb{R}^n$ (when $\kappa =0$) or the Poincar\'e 
  hyperbolic ball ${D}^n$ (when $\kappa < 0$).  In both cases the covering map 
  is a local isometry and $W$ lifts to a infinitesimal homothety $\tilde W$ of 
  $(\tilde M, \tilde h)$.  Since $|\tilde W| < 1$, $(\tilde h , \tilde W)$ 
  determines a globally defined Randers metric $\tilde F$ on $\tilde M$.  
  The projection map is a local isometry of $(\tilde M , \tilde F) \to (M,F)$;
  see Lemma 2 of \cite{BRS04}.
 
  The only infinitesimal homotheties of Euclidean space satisfying 
  $|\tilde W|<1$ globally are the translations.  Hence 
  $(\mathbb{R}^n, \tilde F$), and therefore $(M,F)$, is locally Minkowski.
  (See \cite{BCS00} for the definition of locally Minkowski metrics, and 
   \cite{BRS04} for examples of locally Minkowski Randers metrics.)

  If $\kappa < 0$, then the universal cover is the Poincar\'e ball and
  $\sigma$ must vanish.  So $W$ and $\tilde W$ are infinitesimal isometries.  
  The length of any nonzero, globally defined infinitesimal isometry is 
  unbounded on the ball.  Whence the condition $|\tilde W| = |W| < 1$ 
  forces $\tilde W = 0$.  As a consequence $W = 0$, and $F = h$ is Riemannian.
  \end{proof}


\section{$S^n$ and the Katok examples}
\label{sec:sphere}

According to Theorem \ref{thm:class} the Randers metrics of constant flag 
curvature $K=1$ on the sphere arise as solutions to Zermelo's problem of 
navigation on the canonical Riemannian sphere $(S^n , h)$, under the 
influence of an infinitesimal isometry $W$.  These metrics are projectively 
flat only in the case that $W$ is identically zero \cite{BRS04}; that is, 
$F = h$ is Riemannian.  (Projectively flat Finsler metrics of constant flag 
curvature on the 2-sphere are analyzed in \cite{Br97}.  That paper
is a sequel to \cite{Br96}, Bryant's analysis of arbitrary (not necessarily 
Randers) constant flag curvature metrics on $S^2$.)

Modulo an orthogonal transformation of $\mathbb{R}^{n+1}$, the vector field 
$W$ is expressed by matrix multiplication $W(p) = p \, \Omega$; 
where $p = (p^0 , \ldots , p^n ) \in S^n \subset \mathbb{R}^{n+1}$, and 
$\Omega$ is a $(n+1) \times (n+1)$ constant matrix of the form
\begin{center}
\begin{tabular}{ll}
  $\Omega = a_1 J \oplus \cdots \oplus a_m J \oplus 0$
  & when $n = 2m$ is even,  \\
  $\Omega = a_1 J \oplus \cdots \oplus a_m J$
  & when $n = 2m - 1$ is odd.
\end{tabular}
\end{center}
The $a_i$ are ordered, $a_1 \ge \cdots \ge a_m \ge 0$, and $J$ denotes the 
$2 \times 2$ matrix
\begin{displaymath}
  J = 
  \left(
  \begin{array}{rc}
    0 & 1  \\
   -1 & 0
  \end{array}
  \right) \, .
\end{displaymath}
Note that the globally defined $W$ will have norm less than 1 if and only if 
$a_1 < 1$.  Let's assume this is the case, so that $M = S^n$.

The flow of $W$ is similarly given 
by $\varphi(t,p) = p \, \hbox{Rot}( a_1 t , \ldots , a_m t)$, where 
$\hbox{Rot}(\cdots)$ is the block diagonal matrix
\begin{center}
\begin{tabular}{ll}
  $\hbox{Rot}( a_1 t , \ldots , a_m t) 
   = R(a_1 t) \oplus \cdots \oplus R(a_m t) \oplus 1$,
  & when $n = 2m$ is even;  \\
  $\hbox{Rot}( a_1 t , \ldots , a_m t) 
   = R(a_1 t) \oplus \cdots \oplus R(a_m t)$,
  & when $n = 2m - 1$ is odd.
\end{tabular}
\end{center}
Here, 
\begin{displaymath}
  R( a_i t ) = 
  \left(
  \begin{array}{rc}
    \cos a_i t & \sin a_i t  \\
   -\sin a_i t & \cos a_i t
  \end{array}
  \right) \, .
\end{displaymath}

The Randers metrics generated by the Euclidean $(S^n,h)$ and the 
infinitesimal isometries $W$ via Zermelo navigation were initially introduced 
by Katok \cite{K73}, and later studied by Ziller \cite{Zi82}, in the context 
of Hamiltonian systems.  With a little thought, the following three
observations may be made.

First, if each $a_i$ is rational, then all geodesics of the Randers metric 
close.

Second, given one such Randers metric $F$, Theorem \ref{thm:class} implies it 
is of constant flag curvature $K=1$.  Hence, given any fixed point 
$p \in S^n$, the $F$-distance from $p$ satisfies $d_F(p,q) \le \pi$.  (This is 
the Bonnet-Myers theorem for Finsler metrics.  See \cite{A55,BCS00}.)
Moreover, equality holds for a unique $q$ (cf. Theorem 0.1 of \cite{S96}). 
In the case that $F$ is Riemannian (equivalently, $W=0$) this unique $q$ is 
the antipodal point $-p$.  In general, $q = \varphi(\pi,-p)$.  In analogy with
the Riemannian case, $q$ is the unique point conjugate to $p$ with respect to 
the metric $F$ (Section \ref{sec:conjugate_points}).
  
Last, denote by $\mathcal S^h_r(p) = \{ q \in S^n \ | \ d_h(p,q) = r \}$ the 
{\it geodesic sphere} of radius $r$ about $p$.  Then $\mathcal S^h_r(p)$ 
is a Euclidean $(n-1)$-sphere of radius $\sin(r)$.  The Randers geodesic 
sphere is simply a rotation of $\mathcal S^h_r(p)$.  Explicitly, 
$\mathcal S^F_r(p) = \varphi( r , \mathcal S^h_r(p) )$.   

Suppose $\rho(t)$ is a geodesic of $S^n$, $n = 2m$ or $n=2m-1$, 
invariant under the flow $\varphi(t,\cdot)$.  These include the $m$ geodesics 
parameterizing the great circles 
$C_i = \{ p \in S^n \ | 
        \ p^j = 0\, , \ \forall \ j \not= 2i-1 , 2i \}$, $i = 1, \ldots, m$.
The resulting $F$-geodesic $\P$ is simply a reparameterization of $C_i$.  When 
$\rho(t)$ parameterizes $C_i$ in the direction of the rotation $R(a_i t)$, 
$\P$ has $F$-length $2 \pi / (1 + a_i)$.  Otherwise $\P$ has length 
$2 \pi / (1 - a_i)$.  These two geodesics are considered distinct.  In 
particular, there exist at least $2m$ closed geodesics of $F$ on 
$S^n$, $n=2m-1$, $2m$.   

These geodesics are simple.  In general, $(S^n,F)$ will have geodesics that 
self-intersect.  (This assumes of course that $F$ is non-Riemannian: 
$W \not=0$.)  For example, suppose $p$ is a fixed point of the flow $\varphi$.
In the case that $n=2m$ is even, the north and south poles 
$(0 , \ldots , 0 , \pm 1)$ are always fixed points.  Let 
$\rho : ( -\infty , \infty ) \to S^{2m}$ be the unique $h$-geodesic passing 
through the north pole,
$\rho(2\pi\ell) = ( 0 , \ldots , 0, 1)$, $\ell \in \mathbb{Z}$, with
$\drho(2\pi\ell) = (1,0,\ldots,0)$.  The resulting $F$-geodesic
\begin{eqnarray*}
  \P(t) & = & \varphi( t , \rho(t) )            
        \ = \ ( \sin t , 0 , \ldots , 0 , \cos t ) \, 
              \hbox{Rot}(a_1 t , \ldots , a_m t)        \\
        & = & ( \sin t \cos a_1 t , \sin t \sin a_1 t , 
                0 , \ldots , 0 , \cos t )
\end{eqnarray*}
also passes through the north pole at $t = 2 \pi \ell$, $\ell \in \mathbb{Z}$.
However, 
\begin{center}
  $\dP(2\pi\ell) = ( \cos 2\pi\ell a_1 , \sin 2\pi\ell a_1 , 0 , \ldots , 0 )$.
\end{center}
So $\dP(2\pi\ell) = \dP(0)$ if and only if $a_1 \ell \in \mathbb{Z}$.
Since $0 < a_1 < 1$, $\dP(2\pi) \not= \dP(0)$.  Hence $\P$ self-intersects, 
but does not close, at $t = 2 \pi$.  

Implicit in this discussion is the fact that this $\P$ will close if and only 
if $\ell a_1 \in \mathbb{Z}$ for some integer $\ell$.  That is, $a_1$ must be 
rational.  Indeed, Ziller observed that  
    the $a_i$ may be selected so that the only closed geodesics of $F$ are, 
    modulo parameterization, the $m$ great circles $C_i$.  (As above $n = 2m-1$,
    $2m$.)  The following choice of $a_i$ effects this scenario.  
    Set $a_i = a / p_i$, where $0 < a < 1$ is irrational and the $p_i$ are 
    relatively prime integers with $1 < p_1 < p_2 < \cdots < p_m$.
    Note that there are infinitely many possible choices for the $m$-tuples 
    $a = \{a_i\}_{i=1}^m$, each determining a Randers metric $F_a$.
    Proposition 8 of \cite{BRS04} implies 
    $F_a$ and $F_{\tilde a}$ are locally isometric if and only if
    $a = \tilde a$.  

    Recollect that the $F$-length of $C_i$ depends on the choice of 
    orientation.  We say each orientation produces a distinct geodesic.  
    Hence $(S^n , F)$ has precisely $2m$ closed geodesics.  As a consequence 
    we have the following

\begin{proposition}[\cite{Zi82}]
\label{prop:finite_geos}
  There exist infinitely many, non-isometric Randers metrics of 
  constant flag curvature $K=1$ on $S^{n}$ with only $2m$ 
  closed geodesics.  (Here $n = 2m$ or $2m-1$.)  
\end{proposition}

I would like to make two comments here.  First, 
    the fact that Randers metrics are not reversible is key.  R. Bryant has 
    shown 
    that any reversible Finsler metric of constant flag curvature $K = 1$
    on the 2-sphere is necessarily the standard Riemannian metric of constant 
    sectional curvature 1 \cite{Br04}.  In particular, {\it all} the 
    geodesics close.

    Second, the proposition
    implies there exists a Finsler metric on $S^2$ with only 2 closed 
    geodesics.  This is a {\it non-Riemannian} counter-example to the Three 
    Closed Geodesic Theorem of Lusternik and Schnirelmann \cite{LS29,LS30}: 
    {\it On $S^2$ with an arbitrary Riemannian metric, there exist three 
         simple closed geodesics.}  (See also \cite{B78, G89}.)

The phenomenon of Proposition \ref{prop:finite_geos} is most easily 
illustrated on the 2-sphere, where 
$W(p) = a ( p^1 , - p^0 , p^2 )$, $0 < a < 1$, is an infinitesimal rotation 
about the $p^2$-axis.  The unique unit-speed geodesic of $F$ through 
$q\in S^2$ in the direction $v \in T_q S^2$ is given by 
\begin{displaymath}
  \P(t) = \varphi(t,\rho(t)) = 
  \left(
    \begin{array}{c}
      \rho^0(t) \cos a t \ - \ \rho^1(t) \sin a t \\
      \rho^0(t) \sin a t \ + \ \rho^1(t) \cos a t \\
      \rho^2(t)
    \end{array}
  \right) \, ,
\end{displaymath}
where
$\rho(t) = ( \rho^0(t) , \rho^1(t) , \rho^2(t) ) 
         = \cos(t) \, q + \sin(t) \, (v-W_q)$. 
In the event that $a$ is irrational, it is not difficult to check that $\P(t)$ 
closes if and only if $\rho(t)$ parameterizes the equator $\{ p^2 = 0 \}$.  
As Ziller suggests, this great circle is the only geodesic of $h$ invariant 
under the flow $\varphi(t,\cdot)$.  The $F$-geodesic $\P$ also traces the 
equator.  In the event that $\rho(t)$ parameterizes the equator in the 
direction of rotation by $\varphi$, then the $F$-length of $\P$ is 
$2 \pi / (1 + a)$.  Otherwise, $\P$ has length $2 \pi / (1-a)$.  These two 
parameterizations of the equator are the only closed geodesics on $(S^2,F)$.

Here are two examples of Randers geodesics on the 2-sphere.  The pictures on 
the left are side-views of the sphere, on level with the equator.  On the right
we look down on the North pole.
\begin{center}
Example 11: $a = \fourth$. \\
{\it
    The Riemannian geodesic (here seen as the vertical line bisecting the 
    sphere) passes through the North and South poles, the fixed points of the 
    flow. The parameterization $\rho$ traces the geodesic curve four times
    before the associated Randers geodesic closes.}  \\
\begin{tabular}{cp{10mm}c}  
    \includegraphics{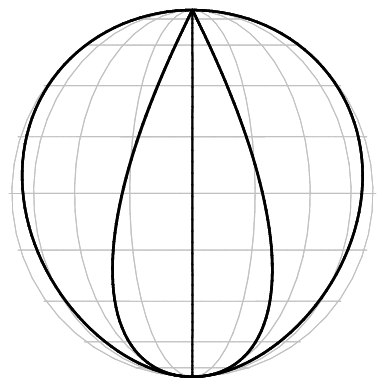} 
  & &
    \includegraphics{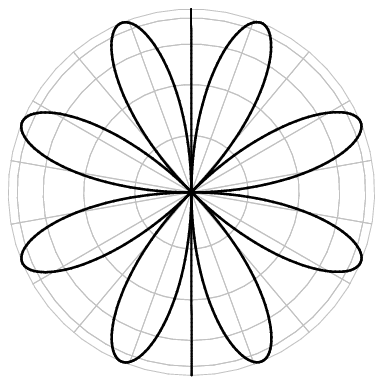} 
\end{tabular}\\
Example 12: $a = \frac{5}{7}$. \\
{\it 
    A Riemannian geodesic (the diagonal line on the left, and ellipse on the
    right) omitting the fixed points of the flow.  The parameterization $\rho$
    circles the sphere seven times before the Randers geodesic closes.} \\
\begin{tabular}{cp{10mm}c}  
    \includegraphics{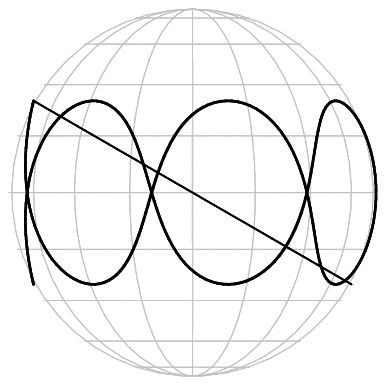}
  & &
  \includegraphics{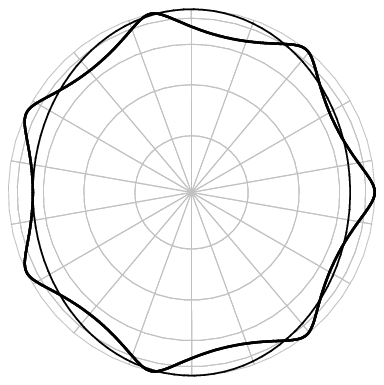}
\end{tabular}
\end{center}

I will close this discussion with two additional observations of Ziller 
\cite{Zi82} and a new result of V. Bangert and Y. Long \cite{BL04}:
    (1)  
    Ziller has shown that any Finsler metric on $S^n$ sufficiently 
    $\mathscr{C}^2$ close to $h$ has {\it at least} $n$ closed geodesics.  
    The Randers examples above show that this lower bound is sharp in the case 
    that $n = 2m$. Bangert--Long have established a stronger result in the 
    case $n=2$: {\it any} Finsler metric on the 2-sphere necessarily has two 
    prime closed geodesics.
    (2)  Ziller also points out that given a generic $\mathscr{C}^2$ Finsler 
    metric on a compact manifold, the initial vectors of closed geodesics are 
    dense in the unit tangent bundle.  In particular, examples of Finsler 
    metrics on $S^n$ with only finitely many closed geodesics are rare.    

It should be noted that Randers metrics of constant curvature on the 
sphere have been studied from several distinct perspectives.  These include
the Bao--Shen \cite{BS02} one-parameter family of examples developed from the 
Hopf fibration on $S^3$, and the approach of Bejancu--Farran through Sasakian 
space forms \cite{BF02, BF03}.

The final two examples in this section are generated from the same Riemannian 
geodesic, but with opposite orientations.  Again the picture on the left is 
a side-view of the sphere, on level with the equator; on the right we look
down on the North pole.
\begin{center}
Example 13: $a = \frac{5}{6}$.  \\
{\it 
    The Riemannian geodesic $\rho$ (seen as a diagonal line on the left, and 
    ellipse on the right) is traced six times before the Randers geodesic 
    $\P_+(t) = \varphi(t,\rho(t))$ closes.  Notice that $\P_+$ is a simple 
    closed curve.} \\  
\begin{tabular}{cp{10mm}c}  
    \includegraphics{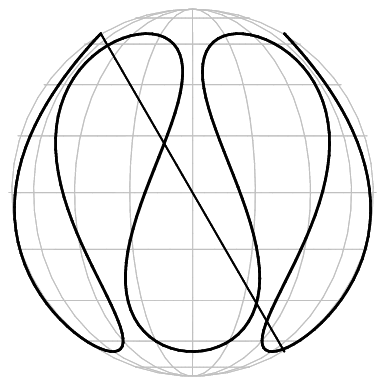} 
  & &
    \includegraphics{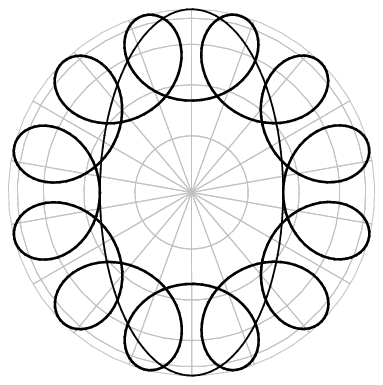}
\end{tabular}\\
Example 14: $a = \frac{5}{6}$.\\
{\it
    This Randers geodesic is generated by tracing the Riemannian geodesic 
    in the opposite direction.  That is, $\P_-(t) = \varphi(t,\rho(-t))$.
    The two Randers geodesics in Examples 13 and 14 are tangent at their
    initial point $\rho(0)$.}\\
\begin{tabular}{cp{10mm}c}  
    \includegraphics{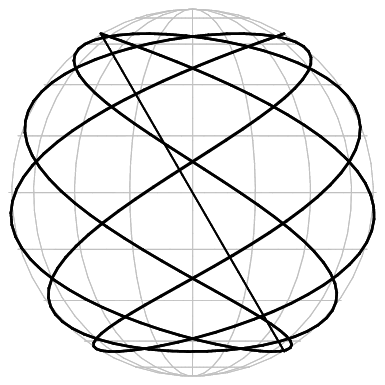} 
  & &
    \includegraphics{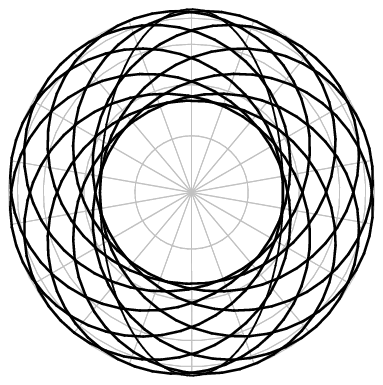}
\end{tabular}
\end{center}

\subsection{A solution to the navigation problem}  I would like to sketch the 
means by which a pilot flying under windy conditions may determine the 
time-efficient paths.  Let me emphasize that {\it this discussion applies only 
to those $W$ generating rigid motions of the sphere.}  Imagine the plane 
traverses the globe (represented by the 2-sphere) with unit speed.  Introduce 
a mild wind (speed less than one) describing rotation about the $z$-axis with 
velocity vector field $W$.  This wind is an isometry, and $|W|<1$, so the 
hypotheses of Theorem \ref{thm:geodesics} are  met.  The plane now travels 
with unit speed relative to the wind, not the earth.  The goal is to determine 
a time-efficient flight path from $p$ to $q$.  

According to the geodesic classification these paths are of the form 
$\P(t) = \varphi(t,\rho(t))$.  So, if $\P : [0,L] \to S^2$ joins $p$ to $q$, 
$\rho(t) = \varphi(-t,\P(t))$ is a geodesic path from $p$ to 
$\varphi(-L,q)$, a point in the pre-orbit of $q$ under the wind $\varphi_t$.  
Note that $L$ is the time it takes the wind to blow a particle from $\rho(L)$ 
to $q$.  These observations suggest the following approach.

Define $q(\tau) = \varphi(-\tau,q)$.  This curve traces the circle in $S^2$ 
parallel to the $xy$-plane, and passing through $q$.  Intuitively, $\tau$ is 
the time it takes the wind to travel from $q(\tau)$ to $q$; and the pilot 
gauges the wind to determine a point $q(\tau)$ at distance $\tau$ from $p$.
Such a point in the pre-orbit certainly exists:  Let $0 \le L(\tau) \le \pi$ 
be the Riemannian distance from $p$ to $q(\tau)$.  If $p \not= q$, then 
$L(0) > 0$ and the intermediate value theorem gives us $\tau_o \in (0,\pi]$ 
such that $L(\tau_o) = \tau_o$.  That is, the Riemannian distance from $p$ to 
$q(\tau_o)$ is equal to the time it takes the wind to blow $q(\tau_o)$ to $q$.
Take the smallest such $\tau_o$.   

Then, if $\rho:[0,L(\tau_o)] \to S^2$ is a unit speed geodesic joining $p$ to
$q(\tau_o)$, $\P(t) = \varphi(t,\rho(t))$ will be a time-minimizing flight 
path from $p$ to $q$.  Note that this scenario is akin to the swimmer at the 
edge of a river, wishing to reach the bank directly opposite, who aims not 
for her desired destination but slightly upstream.

Here is a sequence of pictures illustrating this construction for three 
values of $\tau$: $0 < \tau_- < \tau_o < \tau_+$.  Each image contains the 
reverse flow $\varphi(-t,q)$ from $q$ to $q(\tau)$, the Riemannian geodesic 
$\rho(t) = \rho_\tau(t)$ from $p$ toward $q(\tau) = \varphi(-\tau,q)$, and the 
resulting Randers geodesic $\P(t) = \varphi(t,\rho(t))$; $t \in [0,\tau]$.
\\
\begin{tabular}{ccc}
  \includegraphics{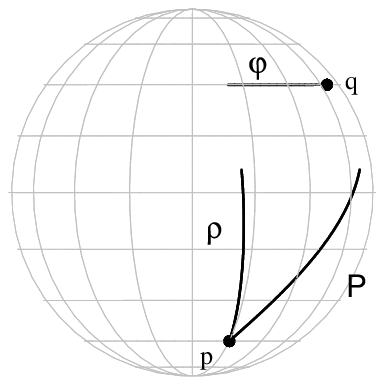} &
  \includegraphics{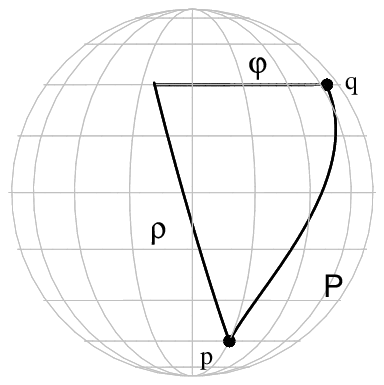} &
  \includegraphics{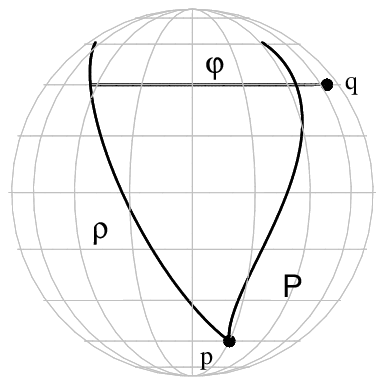} 
\end{tabular}


\section{Conjugate points}
\label{sec:conjugate_points}

Fix $p \in M$.  This section develops the correspondence between points 
$\hat q \in \mathcal M$ conjugate to $p$ with respect to the Riemannian metric 
$h$, and points $q \in M$ conjugate to $p$ with respect to the Randers metric 
$F$.  As before $\P : [0,\ell] \to M$ denotes a unit speed $F$ geodesic 
joining $p$ to $q$. For convenience, I will assume throughout that 
$\rho(t) = \varphi(-t,\P(t))$ is a well-defined curve mapping the interval 
into $\mathcal M$.  (This is the case for all the examples in Sections 
\ref{sec:pictures} and \ref{sec:sphere}.) 

Begin by supposing that $\rho(\ell) = \hat q$ is conjugate to $p$ along $\rho$.
Let $\hat \rho : [ 0 , \hat \ell ] \to \mathcal M$ be a unit $h$-speed
reparameterization of $\rho$.  By assumption there exists a nonzero Jacobi 
field $\mathscr J$ along $\hat \rho$ that vanishes at $p$ and $\hat q$.  
Suppose $\hat \rho_s : [ 0 , \hat \ell ] \to \mathcal M$, 
$-\varepsilon < s < \varepsilon$, is a geodesic variation of 
$\hat \rho = \hat \rho_0$ with variational vector field 
\begin{displaymath}
  \left. \frac{\pa \hat \rho_s}{\pa s} \right|_{s=0} = \mathscr J \, .
\end{displaymath}

I wish to use this Riemannian geodesic variation to construct a Randers 
geodesic variation.  To that end, it is first necessary to reparameterize 
the $\hat \rho_s$ so that the composition with the flow $\varphi$ is a Randers 
geodesic.  Define $\tau : [ 0 , \ell ] \to [ 0 , \hat \ell ]$ by 
\begin{displaymath}
  \tau(t) = \left\{
            \begin{array}{ccl}
            \frac{2}{\sigma} 
            \left( 
              1 - e^{-\sigma t / 2}
            \right) \, , & & \sigma \not= 0 \\
            t \, , & & \sigma = 0 \, . 
            \end{array}
            \right.
\end{displaymath}
Then $\rho_s(t) = \hat \rho_s(\tau(t))$ reparameterizes the geodesics with 
speed $|\drho_s(t)| = v(s) e^{-\sigma t/2}$; here $v(s)$ denotes the constant
$h$-speed of $\hat \rho_s$.

Now $\P_s(t) = \varphi( v(s)t , \rho_s(t) )$ is an $F$-geodesic variation.  
Each geodesic $\P_s$ has constant $F$-speed $v(s)$.  In particular, 
$\P = \P_0$ has speed $v(0) = 1$.  It follows that 
\begin{displaymath}
  J := \left. \frac{\pa \P_s}{\pa s} \right|_{s=0}
     = v'(0) W_\P + (d \varphi) \mathscr J \, 
\end{displaymath}
is a Jacobi field of $F$.  (See \cite{BCS00} for a through treatment of Jacobi 
fields in Finsler geometry).  Making use of the facts that 
$\mathscr J(0) = \mathscr J_p = 0$ and $\mathscr J(\ell) = \mathscr J_{\hat q} = 0$, a standard Jacobi field 
argument yields $v'(0) = 0$.  Hence $J = (d \varphi) \mathscr J$.  Recollect 
that $d \varphi$ denotes the differential of 
$\varphi(t,\cdot)$.  Since $\varphi(t, \cdot)$ is a local 
diffeomorphism, the Jacobian $d \varphi$ is nonsingular.  Hence 
$J$ vanishes if and only if and only if $\mathscr J$ does.  In particular, $J$ 
is nonzero, and $J(0) = 0 = J (\ell)$.  Equivalently, $q = \P(\ell)$
is conjugate to $p = \P(0)$ along $\P$ with respect to the metric $F$.

This establishes the following:  Suppose $\P : [ 0 , \ell ] \to M$ is a 
unit-speed $F$ geodesic with associated Riemannian geodesic $\rho$ (see 
Theorem \ref{thm:geodesics}).  Then $\P(\ell)$ is conjugate to $\P(0)$ 
with respect to $F$ whenever $\rho(\ell)$ is conjugate to 
$\rho(0)$ with respect to $h$.  A similar argument shows the converse holds 
as well, and we have the following
\begin{proposition}
\label{prop:conjugate_pts}
  Assume $(\mathcal M,h)$ is a Riemannian manifold equipped with a globally 
  defined infinitesimal homothety $W$.  Let $F$ denote the Randers metric on 
  $M = \{ |W| < 1 \} \subset \mathcal M$ with navigation data $(h,W)$.  
  Suppose that $\P : [ 0 , \ell ] \to M$ is a geodesic of $(M,F)$, and that 
  the associated Riemannian geodesic $\rho(t) = \varphi(-t , \P(t))$ (Theorem 
  \ref{thm:geodesics}) is defined for all 
  $t \in [0,\ell]$.  Then $\P(\ell)$ is conjugate to $\P(0)$ along $\P$ (with 
  respect to the $F$) if and only if $\rho(\ell)$ is conjugate to $\rho(0)$ 
  (with respect to $h$).
\end{proposition}
\noindent

 In Riemannian geometry $p$ is conjugate to $\hat q$ if and only if $\hat q$ 
 is conjugate to $p$.  This is because $\rho(t)$ is a geodesic if and only if 
 $\rho(-t)$ is as well.  As indicated in the discussion of Case I in Subsection
 $\ref{sec:plane}$, the phenomenon does not hold for Randers metrics in 
 general.  


\section{Globally minimal geodesics}
\label{sec:global_minimizers}

In this section I will establish a relationship between the minimal geodesics 
of $F$ and those of $h$.  Assume throughout that $(\mathcal M, h)$ is complete;
so that $\rho(t)$ and $\varphi(t,p)$ are defined for all $t \in \mathbb{R}$ and
$p \in \mathcal M$.  All the examples of Sections \ref{sec:pictures} and 
\ref{sec:sphere} meet this criterion.  The Hopf--Rinow Theorem assures us that 
any two points of $\mathcal M$ may be joined by a globally $h$-length 
minimizing geodesic.  We will see that, in the event that $M$ satisfies a 
given condition (determined by $\sigma$ and specified below), any two points 
of $(M,F)$ can also be joined by a globally $F$-length minimizing geodesic.

Notice that in each of the examples in \S3 and 4 $\{ |W| < 1 \}$ is connected, 
and {\it convex} in the following sense: any globally $h$-length minimizing 
geodesic joining $p,q \in \{ |W| < 1 \}$, {\it a priori} contained in 
$\mathcal M$, lies in $\{ |W| < 1 \}$.  Moreover, in the example with nonzero 
$\sigma$ (Case II, page \pageref{CaseII}), 
$d_s{}^2 | W_{\eta(s)} |^2 \ge \half \sigma^2$, for any unit speed geodesic 
$\eta(s)$ of $h$.  These observations motivate the hypotheses below.
\begin{proposition}
\label{prop:gm_exist}
Suppose $(\mathcal M,h)$ is a complete, connected Riemannian manifold, and 
that $F$ solves the navigation problem for an infinitesimal homothety $W$ on a 
connected component $M$ of the open sub-manifold 
$\{ |W| < 1 \} \subset \mathcal M$.
\begin{itemize}
  \item[$\circ$]  If $\sigma = 0$, assume $M$ is convex.
  \item[$\circ$]  If $\sigma \not= 0$, assume 
         \begin{displaymath}
           \frac{d^2}{ds^2} |W_{\eta(s)}|^2 \ge \frac{1}{2} \sigma^2 >0
         \end{displaymath}
         for every unit speed $h$-geodesic $\eta$.
\end{itemize}  
Then any two points $p,q$ in the Randers space $(M,F)$ may be joined 
by a globally length minimizing geodesic.
\end{proposition}

{\it Remark:}  Note that every constant flag curvature Randers metric $F$ 
meets the conditions of the proposition.  This is a consequence of 
Theorem \ref{thm:class} which assures us that these metrics
arise as solutions to Zermelo's problem on a Riemannian manifold 
$(\mathcal M ,h)$ of constant curvature under an infinitesimal homothety.  
In the case that $\sigma = 0$ (so that $W$ is an infinitesimal 
isometry) the connected components of $\{ |W| < 1 \}$ are convex in each of 
the three models of Riemannian space forms: the sphere, Euclidean space and 
the hyperbolic ball.  In general, $(\mathcal M,h)$ will be locally isometric 
to one of these three manifolds so that every $p \in M$ admits a convex 
neighborhood.  And so the proposition will hold locally, at the very least,
for these $F$, and globally when $\mathcal M$ is one of the three models.

If, on the other hand $\sigma \not= 0$, then $h$ must be flat and there 
exist coordinates in which $h$ is the standard Euclidean metric and 
$W = \half \sigma x + Qx + C$, for some skew-symmetric matrix $Q$ and 
$C \in \mathbb{R}^n$ (cf. \cite{BRS04}).  Hence 
$d_s{}^2 |W|^2 \ge \half \sigma^2$, and the conditions of the proposition are 
again met.  I will discuss at the end of this section the generality with 
which we may expect this inequality to hold.

\begin{proof}
Let $\ell = d_F(p,q)$ denote the distance from $p$ to $q$ in $M$.
Then there are $F$-unit speed curves $C_i : [ 0 , \ell_i ] \to M$ such that 
$C_i(0) = p$, $C_i(\ell_i) = q$ and $\ell_i \to \ell$.  Define 
$c_i : [0,\ell_i] \to \mathcal{M}$ by $c_i(t) = \varphi(-t,C_i(t))$.  Note 
that 
\begin{displaymath}
  c_i(\ell_i) = \varphi( -\ell_i , C_i(\ell_i) ) = \varphi( -\ell_i , q )
  \longrightarrow \varphi(-\ell,q) \, .
\end{displaymath}
Additionally, Lemma \ref{lem:length}, $F(\dot C_i) = 1$ and the fact that 
$\varphi(t,\cdot)$ is a homothety imply 
$|\dot c_i(t)|^2 = e^{-\sigma t}$.  As a result 
\begin{displaymath}
  d_h(p, c_i(\ell_i)) \le 
  h\hbox{-length}(c_i) = \left\{
                    \begin{array}{ccl}
                      \frac{2}{\sigma} 
                      \left( 1 - e^{-\sigma \ell_i/2} \right) \, ,
                        & & \sigma \not= 0 \\
                      \ell_i \, , 
                        & & \sigma = 0 \, .
                    \end{array}
                      \right.
\end{displaymath}
The continuity of $x \mapsto d_h(p,x)$, in conjunction with these two 
observations, allows us to deduce  
\begin{displaymath}
  d_h(p , \varphi( -\ell , q ) ) \le 
                  \left\{
                    \begin{array}{ccl}
                      \frac{2}{\sigma} 
                       \left( 1 - e^{-\sigma \ell/2} \right) \, ,
                        & & \sigma \not= 0 \\
                      \ell \, ,
                        & & \sigma = 0 \, .
                    \end{array}
                  \right.
\end{displaymath}

Define $q(\tau) := \varphi(-\tau,q)$.  Intuitively, $\tau$ is the time it 
takes the flow to carry $q(\tau)$ onto $q$.  Note that $q = q(0)$.
Since $(\mathcal{M},h)$ is complete, the Hopf-Rinow theorem guarantees a 
globally $h$-length minimizing geodesic 
$\rho_\tau : [0,L(\tau)] \to \mathcal{M}$ joining $p$ to $q(\tau)$, and 
parameterized so that $|\drho_\tau(t)|^2 = e^{-\sigma t}$.  The expression 
above for $d_h(p, \varphi(-\ell,q)) = d_h(p,q(\ell))$ implies 
$L(\ell) \le \ell$.

If $p = q$, the proposition is immediate.  So suppose that 
$p \not= q = q(0)$.  
Then $L(0) > 0$.  The Mean Value Theorem gives us $\ell_o \in (0,\ell]$ 
such that $L(\ell_o) = \ell_o$.  Write $\rho_{\ell_o} = \rho$.  It follows 
that $\P(t) = \varphi(t,\rho(t))$ is a unit speed $F$-geodesic 
joining $\P(0) = p$ to $\P(\ell_o) = \varphi(\ell_o , \rho(\ell_o) ) =
\varphi(\ell_o , q(\ell_o)) = q$.  

Since $\P$ {\it a priori} maps $[0,\ell_o]$ into 
$\mathcal{M}$, it remains to confirm that the interval is mapped into 
$M$.  To that end, consider the sets 
\begin{displaymath}
  M_t := \varphi(t,M) = 
  \{ x \in \mathcal M \ : \ |W_x| < e^{\sigma t/2} \} \, .
\end{displaymath}
The second equality follows from the hypotheses that (a) $\mathcal M$ is 
connected; (b) $d \varphi W_x = W_{\varphi(t,x)}$ (since $\varphi$ is the flow
of $W$); and (c) $W$ is an infinitesimal homothety, 
$|W_{\varphi(t,x)}| = |d \varphi W_x| = e^{\sigma t/2} |W_x|$.  Recall 
$d \varphi$ is the differential of 
$\varphi(t,\cdot): \mathcal M \to \mathcal M$.  In general,
$\P(t) = \varphi(t,\rho(t))$ will lie in $M$ if and 
only if $\rho(t) \in M_{-t}$.   

In the case that $\sigma = 0$, $M_t = M$.  In particular, 
$\rho(0), \rho(\ell_o) \in M = M_{-\ell_o}$, and convexity implies 
$\rho(t) \in M$.  Therefore $\P(t) = \varphi(t,\rho(t)) \in M_{t} = M$.  

Let's now address the case $\sigma \not= 0$.  In order to show that 
$P(t) \in M$, it suffices to see that $\rho(t) \in M_{-t}$.  Equivalently, 
\begin{displaymath}
  |W_{\rho(t)}|^2 < e^{-\sigma t} \, .
\end{displaymath}  
Define $s = (2/\sigma) ( 1 - e^{-\sigma t/2})$.  Then $\eta(s(t)) = \rho(t)$ 
defines a unit speed reparameterization $\eta(s)$ of $\rho(t)$.  The inequality
above is recast as
\begin{displaymath}
  w(s) := | W_{\eta(s)} |^2 
        < \left( 1 - \frac{\sigma}{2} s \right)^2
       =: f(s) \, .
\end{displaymath}
Since $\eta(0) = \rho(0) = p \in M$ and 
$\eta(s(\ell_o)) = \rho(\ell_o)  = \varphi(-\ell_o, q) \in M_{-\ell_o}$, 
the expression holds at the end points.  Additionally, 
$0 < \half \sigma^2 = f''(s) \le w''(s)$.  Therefore the inequality is 
satisfied for all $s \in [0 , s(\ell_o)]$.  We conclude $\P(t) \in M$.  Whence 
$\P$ is a globally length minimizing geodesic of $(M,F)$.  Finally, 
$d_F(p,q) = \ell$ implies $\ell = \ell_o$.
\end{proof}

Implicit in this proof is the fact that if $\P$ is a globally minimal geodesic
of $(M,F)$, then the associated Riemannian geodesic $\rho$ is a global 
minimizer of $(\mathcal M,h)$.  The converse holds as well.

Suppose that $\rho$ is a length minimizing geodesic of $(\mathcal M , h)$ 
joining $\rho(t) = p$ to $\rho(\ell) = \hat q$.  If 
$\P(t) = \varphi(t,\rho(t))$ lies in $M$ for all $t$, so that it is a geodesic 
of $(M,F)$, then $\P$ is also a global minimizer.  To see this, we argue by 
contradiction, assuming that $d_F(p,q) = \ell_o < \ell$, where $q := \P(\ell)$.
Let $\P_o :[0,\ell_o] \to M$ be a length minimizing curve joining 
$p = \P_o(0)$ to $q = \P_o(\ell_o)$, and consider the associated Riemannian 
geodesic $\rho_o(t) = \varphi(-t,\P_o(t))$ from $p$ to $\rho_o(\ell_o)$.  
 The path $\xi : [- \ell_o , -\ell] \to \mathcal M$
defined by $\xi(t) = \varphi(t,q)$ runs from $\rho_o(\ell_o)$ to $\rho(\ell)$.
Consequently $\hbox{length}_h(\rho) \le \hbox{length}_h(\rho_o) + 
\hbox{length}_h(\xi)$.  Computing $|\dot\xi(t)|^2 = |W_{\varphi(t,q)}|^2 = 
|d \varphi W_q|^2 = e^{\sigma t} |W_q|^2 < e^{\sigma t}$, a straightforward 
computation produces
$\hbox{length}_h(\xi) < \hbox{length}_h(\rho) - \hbox{length}_h(\rho_o)$. 
This is a contradiction; consequently $\P$ must be minimal.
We have established the following
\begin{proposition}
\label{prop:global_minimizers}
  Suppose $(\mathcal M,h)$ is a complete, connected Riemannian manifold, and 
  that $F$ solves the navigation problem for an infinitesimal homothety $W$ on 
  a connected component $M$ of the open sub-manifold $\{ |W| < 1 \}
  \subset \mathcal M$.  
  \begin{itemize}
    \item[$\circ$]
      If the homothety constant $\sigma$ is zero, assume $M$ is convex.  
    \item[$\circ$]
      If $\sigma \not= 0$, assume that
         \begin{displaymath}
           \frac{d^2}{ds^2} |W_{\eta(s)}|^2 \ge \frac{1}{2} \sigma^2 >0
         \end{displaymath}
         for every unit speed $h$ geodesic $\eta$.
  \end{itemize}  
  Then a geodesic $\P$ of $F$ is a global minimizer of 
  $F$-path-length from $\P(0)$ to $\P(\ell)$ if and only if the associated 
  Riemannian geodesic $\rho$ (Theorem \ref{thm:geodesics}) globally minimizes 
  $h$-path-length between $\rho(0)$ and $\rho(\ell)$.
\end{proposition}
\begin{corollary}
The point $q$ is a cut point of $p$, with respect to $F$, if and only if 
$\hat q = \varphi(-\ell,q)$ is a cut point of $p$ with respect to $h$.
(Here $\ell$ is the $F$-distance from $p$ to $q$.)
\end{corollary}

As promised, I close this section by discussing the inequality 
$d_s{}^2 |W|^2 \ge \half \sigma^2$.  First note that it does not hold for 
any infinitesimal homothety.  As a counter-example
take the standard 2-sphere parameterized by 
\begin{center}
$(\theta, \phi) \mapsto ( \cos \theta \, \sin \phi , \sin \theta \, \sin \phi ,
 \cos \phi)$.  
\end{center}
Note that $W = \partial_\theta$ is an infinitesimal 
isometry, and $\eta(s) = ( \theta(s) , \phi(s) ) = ( 0 , s )$ a unit 
speed geodesic along which $d_s{}^2 |W|^2 = 2 \, (\sin^2(s) - \cos^2(s) )$ 
fails the inequality.

However, the inequality will hold whenever the Riemannian sectional curvature 
$\kappa(\dot \eta \wedge W )$ of the plane spanned by $\{ \dot \eta , W \}$ 
is nonpositive.  To see this note that 
\begin{displaymath}
  |\nabla_{\dot \eta} W| \ \ge \ h( \dot \eta , \nabla_{\dot \eta} W )
                         \  =  \ \half \sigma \, ,
\end{displaymath}
a consequence of the homothety hypothesis 
$h( v , \nabla_v W ) = \half \sigma |v|^2$, for all $v \in T_x \mathcal M$.
Whence
\begin{equation}
\label{eqn:ineq2}
  d_s{}^2 |W|^2 \ = \ 2 | \nabla_{\dot \eta} W|^2 
                    \, + \, 2 h( W , \nabla_{\dot \eta}{}^2 W )
                \ \ge \ \half \sigma^2 
                    \, + \, 2 h( W , \nabla_{\dot \eta}^2 W ) \, ,
\end{equation}
Since $W$ is an infinitesimal homothety we have 
\begin{equation}
\label{eqn:ineq1}
  W_{i:j:k} + W_{j:i:k} = 0 \, .
\end{equation}  
Consequently, the inner product appearing above may be re-expressed as 
\begin{eqnarray*}
  h( W , \nabla_{\dot \eta}^2 W ) 
  & = & W^i W_{i:j:k} \dot \eta^j \dot \eta^k  
  \ = \ - \dot \eta^i W_{i:j:k} W^j \dot \eta^k \\
  & = & - \kappa( \dot \eta \wedge W ) \, ( |W|^2 - h(\dot \eta , W)^2 ) \, .
\end{eqnarray*}
Here, Equation \ref{eqn:ineq1} provides the second equality; and the Ricci 
identity for $W$, followed by a second application of (\ref{eqn:ineq1}), 
provides the last.  Since $|W|^2 - h(\dot \eta , W)^2 \ge 0$, we have 
$h(W , \nabla_{\dot \eta}^2 W) \ge 0$ whenever 
$\kappa( \dot \eta \wedge W ) \le 0$, and Equation \ref{eqn:ineq2} implies 
\begin{displaymath}
  d_s{}^2 |W|^2 \ge \half \sigma^2 \, .
\end{displaymath}
\vspace{0.05in}

\noindent
{\bf Acknowledgments}.  I thank D. Bao and W. Ziller for many helpful 
comments.



\end{document}